\documentclass[12pt]{amsart}
\usepackage{mathrsfs, amssymb, bm}
\usepackage[all]{xy}

\newif\ifappendix
\appendixfalse

\newcommand{\Pa}[9]{\bibitem{#1} {#2}, \emph{#3}, {#4} \textbf{#5} ({#6}), {#7}--{#8}.}
\newcommand{\ed}{

\Newpage

\Newpage

\ifappendix
\newpage
\centerline{\textsc{Personal appendix: Three fundamental problems}}
\par\smallskip

Following are two suggested extensions to F761, and one suggested extension for F762.
A subgroup $G$ of $\Bgp$ is \emph{Menger-bounded} iff:\\
There is $f\in\NN$ such that
$$(\forall g\in G)(\exists^\oo n)\ |g|\rest [0,n)\le f(n).$$
$G^2$ is Menger-bounded iff:\\
There is $f\in\NN$ such that
$$(\forall g_0,g_1\in G)(\exists^\oo n)(\forall i=0,1)\ |g_i|\rest [0,n)\le f(n).$$

We used a weak but unprovable hypothesis to prove
that there is a group $G\le\Bgp$ such that
$G$ is Menger-bounded, but $G^{2}$ is not Menger-bounded.
Now assume that we are given more freedom. I expect the following
problem of Tka\v{c}enko to have a positive answer.

\bPrb{F761(A)}
Are there, in ZFC, Menger-bounded groups $G,H\le\Bgp$ such that
$G\x H$ is not Menger-bounded.
\ePrb

\bdfn
A subgroup $G$ of $\Bgp$ is \emph{Rothberger-bounded} iff\\
For each increasing $h\in\NN$, there is $\vphi:\N\to\FinSeqs{\bbZ}$
such that:
$$(\forall g\in G)(\exists n)\ g\rest [0,h(n)) = \vphi(n).$$
\edfn

\bPrb{F761(B)}
Does \textsf{CH} imply the existence of a group $G\le\Bgp$ such that $G$ is
Rothberger-bounded but $G^2$ is not Menger-bounded?
\ePrb

\emph{Semifilter-trichotomy} is the hypothesis equivalent to $\fu<\fg$,
which asserts that for each semifilter $\cF$ on $\N$ (i.e.,
$\cF\sbst\roth$ is nonempty, and for all $A,B\sbst\N$, $\cF\ni A\as B\to B\in\cF$),
there is an increasing sequence $h$ such that $\cF/h$ is either the
Fr\'echet filter (all cofinite sets), or an ultrafilter, or $\roth$.

\bPrb{F762(C)}
Does semifilter-trichotomy imply that the square of each
Menger-bounded subgroup of $\Bgp$ is Menger-bounded?
\ePrb

The question for \emph{larger} powers was settled in the positive by
Banakh and Zdomskyy, and independently bey Heike in her work on
F762.

\newpage

\fi

\end{document}}

\newcommand{\bPrb}[1]{\par\medskip\noindent\textbf{Problem #1.}}
\newcommand{\ePrb}{\par\medskip}

\newcommand{\mbx}[1]{{\mbox{#1}}}
\long\def\forget#1\forgotten{}
\newcommand{\FinSeqs}[1]{{{#1}^{<\alephes}}}
\newcommand{\beq}[1]{\begin{equation}\label{#1}}
\newcommand{\eeq}{\end{equation}}
\newcommand{\vzero}{{\vec{0}}}
\newcommand{\NNup}{{\N^{\uparrow\N}}}
\newcommand{\mx}[1]{{\left (\begin{matrix}#1\end{matrix}\right )}}
\newcommand{\intvl}[2]{{[#1(#2),#1(#2\!+\!1))}}
\newcommand{\alephes}{{\aleph_0}}
\newcommand{\sseq}[1]{\{#1 : n\in\N\}}
\newcommand{\oo}{\infty}
\newcommand{\Newpage}{}
\newcommand{\rest}{\upharpoonright}
\newcommand{\Bgp}{{\bbZ^\N}}
\newcommand{\CH}{the Continuum Hypothesis}
\newcommand{\x}{\times}

\newcommand{\vphi}{\varphi}

\newcommand{\Impl}{\Rightarrow}
\newcommand{\seq}[1]{\{#1\}_{n\in\N}}

\newcommand{\spst}{\supseteq}
\newcommand{\sms}{{\ }}
\newcommand{\bo}[1]{{\bm{[}\sms{}#1\sms{}\bm{]}}}

\newcommand{\fb}{\mathfrak{b}}
\newcommand{\fd}{\mathfrak{d}}
\newcommand{\fu}{\mathfrak{u}}
\newcommand{\fg}{\mathfrak{g}}
\newcommand{\NN}{{\N^{\N}}}

\newcommand{\itm}{\item}
\newcommand{\cU}{\mathcal{U}}

\newcommand{\cM}{\mathcal{M}}

\newcommand{\op}{\operatorname}

\newcommand{\dom}{\op{dom}}
\newcommand{\cov}{\op{cov}}
\newcommand{\be}{\begin{enumerate}}
\newcommand{\ee}{\end{enumerate}}

\newcommand{\Union}{\bigcup}
\newcommand{\cF}{\mathcal{F}}
\newcommand{\cP}{\mathcal{P}}
\newcommand{\N}{\mathbb{N}}
\newcommand{\roth}{{[\N]^{\aleph_0}}}

\newcommand{\bbQ}{\mathbb{Q}}
\newcommand{\bbZ}{\mathbb{Z}}
\newcommand{\bbC}{\mathbb{C}}
\newcommand{\sbst}{\subseteq}

\newcommand{\as}{\subseteq^*}

\newcommand{\fr}{\mathfrak{r}}

\newtheorem{thm}{Theorem}
\newcommand{\bthm}{\begin{thm}} \newcommand{\ethm}{\end{thm}}
\newtheorem{prop}[thm]{Proposition}
\newcommand{\bprp}{\begin{prop}} \newcommand{\eprp}{\end{prop}}
\newtheorem{fact}[thm]{Fact}
\newcommand{\bfct}{\begin{fact}} \newcommand{\efct}{\end{fact}}
\newtheorem{prob}[thm]{Problem}
\newcommand{\bprb}{\begin{prob}} \newcommand{\eprb}{\end{prob}}
\newtheorem{lem}[thm]{Lemma}
\newcommand{\blem}{\begin{lem}} \newcommand{\elem}{\end{lem}}
\newtheorem{cor}[thm]{Corollary}
\newcommand{\bcor}{\begin{cor}} \newcommand{\ecor}{\end{cor}}
\newtheorem{conj}[thm]{Conjecture}
\newcommand{\bcnj}{\begin{conj}} \newcommand{\ecnj}{\end{conj}}
\theoremstyle{definition}
\newtheorem{defn}[thm]{Definition}
\newcommand{\bdfn}{\begin{defn}} \newcommand{\edfn}{\end{defn}}
\theoremstyle{remark}
\newtheorem{rem}[thm]{Remark}
\newcommand{\brem}{\begin{rem}} \newcommand{\erem}{\end{rem}}
\newtheorem{exam}[thm]{Example}
\newcommand{\bexs}{\begin{exam}} \newcommand{\eexs}{\end{exam}}
\newcommand{\bpf}{\begin{proof}} \newcommand{\epf}{\end{proof}}

\title{Squares of Menger-bounded groups}

\author{Micha\l{} Machura}
\address{Institute of Mathematics, University of Silesia, ul.\ Bankowa 14, 40-007 Katowice, Poland; and
Departament of Mathematics, Bar-Ilan University, Ramat Gan 52900,
Israel} \email{machura@ux2.math.us.edu.pl}

\author{Saharon Shelah}
\address{Einstein Institute of Mathematics, The Hebrew University of Jerus\-alem,
Givat Ram, 91904 Jerus\-alem, Israel, and Mathematics Department,
Rutgers University, New Brunswick, NJ, USA}
\email{shelah@math.huji.ac.il}

\author{Boaz Tsaban}
\address{Department of Mathematics, Bar-Ilan University, Ramat-Gan 52900, Israel; and
Department of Mathematics, Weizmann Institute of Science, Rehovot 76100, Israel}
\email{tsaban@math.biu.ac.il}
\urladdr{http://www.cs.biu.ac.il/\~{}tsaban}
\thanks{The authors were partially supported by: The EU Research and Training Network HPRN-CT-2002-00287,
United States-Israel BSF Grant 2002323, and
the Koshland Center for Basic Research, respectively.
This is the second author's Publication 903.
}


\begin{document}

\begin{abstract}
Using a portion of \CH{}, we prove that
there is a Menger-bounded (also called $o$-bounded)
subgroup of the Baer-Specker group $\Bgp$, whose square is
not Menger-bounded.
This settles a major open problem concerning boundedness notions
for groups, and implies that Menger-bounded groups need not be Scheepers-bounded.
This also answers some questions of Banakh, Nickolas, and Sanchis.
\end{abstract}

\maketitle

\Newpage

\section{Introduction}

Assume that $(G,\cdot)$ is a topological group.
For $A,B\sbst G$, $A\cdot B$ stands for $\{a\cdot b : a\in A,\ b\in B\}$, and
$a\cdot B$ stands for $\{a\cdot b : b\in B\}$.
The following definitions are due, independently, to Okunev and Ko\v{c}inac.

\begin{defn}\label{bddgps}
Assume that $(G,\cdot)$ is a topological group. $G$ is:
\be

\itm \emph{Menger-bounded} if for each sequence
$\seq{U_n}$ of neighborhoods of the unit, there exist finite sets
$F_n\sbst G$, $n\in\N$, such that $G=\Union_n F_n\cdot U_n$.

\itm \emph{Scheepers-bounded} if for each sequence
$\seq{U_n}$ of neighborhoods of the unit, there exist finite sets
$F_n\sbst G$, $n\in\N$, such that for each finite set $F\sbst G$,
there is $n$ such that $F\sbst F_n\cdot U_n$.

\itm \emph{Hurewicz-bounded} if for each sequence
$\seq{U_n}$ of neighborhoods of the unit, there exist finite sets
$F_n\sbst G$, $n\in\N$, such that for each $g\in G$, $g\in F_n\cdot U_n$
for all but finitely many $n$.

\itm \emph{Rothberger-bounded} if for each sequence
$\seq{U_n}$ of neighborhoods of the unit, there exist elements
$a_n\in G$, $n\in\N$, such that $G=\Union_n a_n\cdot U_n$.

\ee
\end{defn}

Several instances of these properties were studied in, e.g.,
\cite{TkaIntro, Hernandez, HRT, KMexample, BNS, o-bdd}.
A study from a more general point of view was initiated in \cite{Koc03, coc11, Bab05}.
These properties are obtained from the following
general topological properties by
restricting attention to open covers of the form $\{a\cdot U : a\in G\}$,
where $U$ is an open neighborhood of the unit.

\begin{defn}\label{sps}
Assume that $X$ is a topological space. $X$ has the
\be

\itm \emph{Menger property} \cite{Menger24}
if for each sequence
$\seq{\cU_n}$ of open covers of $X$, there exist finite sets
$\cF_n\sbst\cU_n$, $n\in\N$, such that $\Union_{n\in\N} \cF_n$ is a cover of $X$.

\itm \emph{Scheepers property} \cite{coc1}
if for each sequence
$\seq{\cU_n}$ of open covers of $X$, there exist finite sets
$\cF_n\sbst\cU_n$, $n\in\N$, such that
for each finite set $F\sbst X$,
there is $n$ such that $F\sbst \Union_{U\in\cF_n} U$.

\itm \emph{Hurewicz property} \cite{Hure25, Hure27}
if for each sequence
$\seq{\cU_n}$ of open covers of $X$, there exist finite set
$\cF_n\sbst\cU_n$, $n\in\N$, such that
for each element $x\in X$,
$x\in\Union_{U\in\cF_n} U$ for all but finitely many $n$.

\itm \emph{Rothberger property} \cite{Roth38}
if for each sequence
$\seq{\cU_n}$ of open covers of $X$, there exist elements
$U_n\in\cU_n$, $n\in\N$, such that
$X=\Union_{n\in\N} U_n$.
\ee
\end{defn}
Except for the second, all these properties are classical.
They share the same structure and can be defined in a unified manner \cite{coc1, coc2}.
These properties were analyzed in many papers
and form an active area of mathematical research -- see \cite{LecceSurvey, KocSurv, ict, BZSPM}
and references therein.

The relations between the mentioned group theoretic and general topological properties
is thoroughly investigated in \cite{BG, WeissBG}. Here, we consider only the group theoretic
properties.
Clearly, the group theoretic properties are related as follows:
{\small
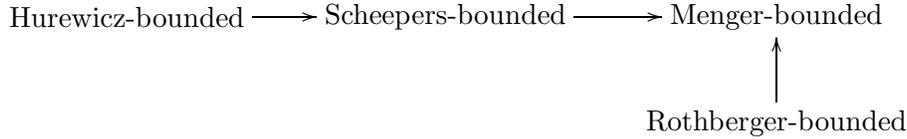
\begin{center}
$\xymatrix{
\mbx{Hurewicz-bounded} \ar[r]        & \mbx{Scheepers-bounded} \ar[r] & \mbx{Menger-bounded}\\
 & & \mbx{Rothberger-bounded}\ar[u]
}$
\end{center}
}
Babinkostova \cite{Bab05} proved that a metrizable group $G$ is Hurewicz-bounded if, and only if,
$G$ is a subgroup of a $\sigma$-compact group (see \cite{BG}).
Neither the leftmost horizontal implication, nor the vertical implication, can be inverted---even when restricting attention to metrizable groups.
The question whether the remaining implication can be inverted
remained thus far open \cite{BNS, coc11, BZSPM, BG}.

\bprb\label{mainprob}
Is every Menger-bounded group Scheepers-boun\-ded?
\eprb

The notions of Menger-bounded and Scheepers-bounded groups are related
in the following elegant manner. For each $k$, let $G^k$ be the direct
product of $k$ copies of $G$.

\bthm[Babinkostova-Ko\v{c}inac-Scheepers \cite{coc11}]\label{powchar}
$G$ is Scheepers-boun\-ded if, and only if, $G^k$ is Menger-bounded
for all $k$.
\ethm

In light of Theorem \ref{powchar}, Problem \ref{mainprob}
asks whether there could be a (metrizable) group $G$ such that for some $k$,
$G^k$ is Menger-bounded but $G^{k+1}$ is not.
The proof of Theorem \ref{powchar} in \cite{coc11} actually shows that the following
holds for each natural number $k$. Since this is used in the sequel, we give a proof.

\blem\label{klem}
$G^k$ is \emph{Menger-bounded} if, and only if,
for each sequence $\seq{U_n}$ of neighborhoods of the unit of $G$,
there exist finite sets
$F_n\sbst G$, $n\in\N$, such that for each $F\sbst G$ with $|F|=k$,
there is $n$ such that $F\sbst F_n\cdot U_n$.
\elem
\bpf
$(\Impl)$ Let $U_n$, $n\in\N$, be neighborhoods of the unit of $G$.
Then $U_n^{\ k}$, $n\in\N$, are neighborhoods of the unit of $G^k$.
Take finite $G_n\sbst G^k$ such that $G^k=\Union_n G_n\cdot U_n^{\ k}$.
Adding elements if necessary, we may assume that each $G_n$ has the form $F_n^{\ k}$
for some finite $F_n\sbst G$. The sets $F_n$ are as required: Given
$g_1,\dots,g_k\in G$, there is $n$ be such that $(g_1,\dots,g_k)\in
F_n^{\ k}\cdot U_n^{\ k} = (F_n\cdot U_n)^k$, and therefore
$g_1,\dots,g_k\in F_n\cdot U_n$.

$(\Leftarrow)$ It suffices to consider basic neighborhoods of the unit of $G^k$.
Let $V_n=U_{n,1}\x \dots\x U_{n,k}$, $n\in\N$, be such that each $U_{n,i}$
is a neighborhood of the unit of $G$. For each $n$, $U_n=U_{n,1}\cap \dots\cap U_{n,k}$
is a neighborhood of the unit of $G$. Take finite $F_n\sbst G$, $n\in\N$, such that
for each $F\sbst G$ with $|F|=k$, there is $n$ such that $F\sbst F_n\cdot V_n$.
Given $(g_1,\dots,g_k)\in G^k$, take $F=\{g_1,\dots,g_k\}$.
If needed, add elements to $F$ to have $|F|=k$.
Then, whenever $F\sbst F_n\cdot U_n$, we have that $(g_1,\dots,g_k)\in F^k\sbst (F_n\cdot U_n)^k\sbst
F_n^{\ k}\cdot V_n$. Thus, the finite sets $F_n^{\ k}$, $n\in\N$, are as required
for the Menger-boundedness of $G^k$.
\epf

We give a negative answer to Problem \ref{mainprob} by showing that,
assuming \CH{} or just a portion of it, there is
\emph{for each} $k$ a metrizable group $G$ such that $G^k$ is Menger-bounded but $G^{k+1}$ is not.

Some special hypothesis is necessary in order to prove such a
result: Banakh and Zdomskyy \cite{BZSPMQdM, BZSPM}, and later (independently)
Mildenberger and Shelah \cite{F761}, proved that
consistently, every topological group with Men\-ger-bounded square
is Scheepers-bounded.

Question 1 of Banakh, Nickolas, and Sanchis \cite{BNS} asks whether each
Menger-bounded subgroup of $\bbC^\N$ (with coordinate-wise addition)
is mixable or $o_\cF$-bounded for some filter $\cF$. As it is proved there that mixable Menger-bounded groups
are Scheepers bounded, and the same holds for groups which are
$o_\cF$-bounded for some filter $\cF$, we obtain a negative answer to both
questions: The groups we construct are, in particular, subgroups of $\bbC^\N$.

\medskip

The problem whether, consistently, every Menger-bounded group is
Scheepers-bounded is yet to be addressed. The answer to this problem is
positive if, and only if, the answer to the following problem is positive.
\bprb\label{stillopen}
Is it consistent that for each Menger-bounded group $G$, $G^2$ is Menger-bounded?
\eprb

There seems to be no straightforward \emph{negative} answer to Problem \ref{stillopen}.
If $G$ abelian and Menger-bounded but $G^2$ is not, then $G$ cannot be analytic,
and not a free topological group over a Tychonoff space, either \cite{BZSPM, ZdFree06}.

\section{Specializing the question for the Baer-Specker group}

The \emph{Baer-Specker group} is the abelian group $(\Bgp,+)$, where $+$ denotes
coordinate-wise addition.
Subgroups of the Baer-Specker group form a rich
source of examples of groups with various boundedness properties \cite{Baer37, Spec50, BlAb, BG, WeissBG}.
The advantage of working in $\Bgp$ is that the boundedness properties
there can be stated in a purely combinatorial manner.

We use mainly self-evident notation.
The quantifiers $(\exists^\oo n)$ and $(\forall^\oo n)$ stand for ``there exist infinitely many $n$'' and
``for all but finitely many $n$'', respectively.
The canonical basis for the topology of $\Bgp$ consists of the sets
$$\bo{s} = \{f\in\Bgp : s\sbst f\}$$
where $s$ ranges over all finite sequences of integers.
For natural numbers $k<m$, $[k,m)=\{k,k+1,\dots,m-1\}$.
For a partial function $f:\N\to\bbZ$, $|f|$ is the function with the same domain,
which satisfies $|f|(n)=|f(n)|$, where in this case $|\cdot|$ denotes the absolute value.
For partial functions $f,g:\N\to\N$ with $\dom(f)\sbst\dom(g)$,
$f\le g$ means: For each $n$ in the domain of $f$,
$f(n)\le g(n)$. Similarly, $f\le k$ means: For each $n$ in the domain of $f$,
$f(n)\le k$. Finally, for a set $X$ and $k\in\N$, $[X]^k=\{F\sbst X: |F|=k\}$.

In a manner similar to the characterizations given in \cite{BG},
we prove the following.

\bthm\label{kMenBdd}
Assume that $G$ is a subgroup of $\Bgp$. The following
conditions are equivalent:
\be
\itm $G^k$ is Menger-bounded.
\itm For each increasing $h\in\NN$, there is $f\in\NN$ such that:
$$(\forall F\in [G]^k)(\exists n)(\forall g\in F)\ |g|\rest [0,h(n))\le f(n).$$
\itm For each increasing $h\in\NN$, there is $f\in\NN$ such that:
$$(\forall F\in [G]^k)(\exists^\oo n)(\forall g\in F)\ |g|\rest [0,h(n))\le f(n).$$
\itm There is $f\in\NN$ such that:
$$(\forall F\in [G]^k)(\exists^\oo n)(\forall g\in F)\ |g|\rest [0,n)\le f(n).$$
\ee
\ethm
\bpf
$(1\Impl 2)$ Fix an increasing $h \in \NN$.
For each $n$, take $U_n = \bo{0\rest [0,h(n))}$.
Using Lemma \ref{klem}, find finite $F_n\sbst G$, $n\in\N$, such that
each $k$-element subset of $G$ is contained in $F_n + U_n$ for some $n$.
Define $f\in\NN$ by
$$f(n) = \max\{|a(i)| : a\in F_n\mbox{ and }i<h(n)\}$$
for each $n$.
Fix $F\in [G]^k$. Take $n$ such that $F\sbst F_n+U_n$.
For each $g\in F$, there is $a\in F_n$ such that $g\in a+U_n = \bo{a\rest[0,h(n))}$,
that is, $g\rest[0,h(n))=a\rest[0,h(n))$, and therefore $|g|\rest [0,h(n)) = |a|\rest[0,h(n))\le f(n)$.

$(2\Impl 1)$
Assume that $\seq{U_n}$ is a sequence of neighborhoods of $0$ in $\Bgp$.
Take an increasing $h\in\NN$ such that $\bo{0\rest[0,h(n))}\sbst U_n$ for each $n$.
Apply (2) for $h$ to obtain $f$.
For each $n$ and each $s\in\bbZ^{[0,h(n))}$ with $|s|\le f$,
choose (if possible) $a_s \in G$ such that
$a_s\rest[0,h(n)) = s$. If this is impossible, take $a_s = 0$.
Let $F_n = \{a_s : s\in\bbZ^{[0,h(n))},\ |s|\le f\}$.
We claim that the sets $F_n$ are as required in Lemma \ref{klem}.
Given $F\in [G]^k$, let $n$ be such that $|g|\rest[0,h(n)) \le f(n)$ for each $g\in F$.
Then for each $g\in F$, there is $s\in\bbZ^{[0,h(n))}$ such that $g\rest[0,h(n))=s=a_s\rest[0,h(n))$,
and thus
$$g\in \bo{a_s\rest[0,h(n))} = a_s+\bo{0\rest[0,h(n))}\sbst a_s+U_n\sbst F_n+U_n.$$

$(1\Impl 3)$ Let $\N=\Union_mI_m$ be a partition into infinite sets.
Fix $m$. For each $n\in I_m$, take $U_n = \bo{0\rest [0,h(n))}$.
The arguments of $(1\Impl 2)$ show that there is $f_m:I_m\to\N$, such that
for each $F\in [G]^k$,
there is $n\in I_m$ such that $|g|\rest [0,h(n)) \le f_m(n)$ for all $g\in F$.
Take $f=\Union_m f_m$.

$(3\Impl 2)$ and $(3\Impl 4)$ are trivial.

$(4\Impl 3)$ This was pointed out by Banakh and Zdomskyy, and later independently
by Simon. Indeed, fix any increasing $h\in\NN$. Let $f$ be as in $(4)$.
We may assume that $f$ is increasing. Define $\tilde f(n)=f(h(n+1))$ for each $n$.
Fix $F\in [G]^k$.
$$I=\{n : n>h(0)\mbox{ and } (\forall g\in F)\ |g|\rest [0,n)\le f(n)\}$$
is infinite.
For each such $n\in I$, let $m$ be such that $n\in\intvl{h}{m}$.
Then for each $g\in F$,
$$|g|\rest [0,h(m))\le |g|\rest [0,n) \le f(n) \le f(h(m+1)) = \tilde f(m).$$
As $I$ is infinite, there are infinitely many such $m$.
\epf

\section{An important corollary of the main theorem}

The purpose of this section is twofold: Making a significant corollary of
our main result (Theorem \ref{main}) accessible to a wider audience,
and exposing the reader to the technically delicate proof of Theorem \ref{main}
via a more accessible proof. Readers who are experienced with cardinal characteristics
of the continuum may, however, wish to try moving directly to the next section, which is
essentially self-contained.

\bthm[$\mathsf{CH}$]\label{mini}
There is a Menger-bounded group $G\le\Bgp$ such that
$G^2$ is not Menger-bounded.
\ethm

We first give an informal outline of the proof.
Assume \CH. By transfinite induction on $\alpha<\aleph_1$,
we will choose generators $g^\alpha_0,g^\alpha_1\in\Bgp$,
and let $G\le\Bgp$ be the group generated by $\{g^\alpha_0,g^\alpha_1 : \alpha<\aleph_1\}$.

Enumerate $\bbZ^2=\{(a_n,b_n): n\in\N\}$ with each pair $(a,b)$ occurring
infinitely often, and enumerate $\Bgp=\{d_\alpha : \alpha<\aleph_1\}$.
At step $\alpha$, let $M_\alpha\le\Bgp$ contain all functions encountered
in earlier steps, as well as $d_\alpha$, and assume that $M_\alpha$ is closed
under all operations required in the proof.
$M_\alpha$ is countable, and we choose $h_\alpha\in\NN$
which grows much faster than any element of $M_\alpha$.
Fix $n$. We choose a solution of $a_nx+b_ny=0$ over $\bbZ$ with $\max\{|x|,|y|\}\ge d_\alpha(h_\alpha(n+1))$,
but not greater than necessary (henceforth: minimal solution).
For each $k\in\intvl{h}{n}$, we set $(g^\alpha_0(k),g^\alpha_1(k))=(x,y)$.

The fact that $\max\{|g^\alpha_0(h_\alpha(n))|,|g^\alpha_1(h_\alpha(n))|\}\ge d_\alpha(h_\alpha(n+1))$
guarantees that $G^2$ is not Menger-bounded (using Theorem \ref{kMenBdd}(4)).

The proof that $G$ is Menger-bounded is more subtle (a preservation argument).
A general element of $G$ is a linear combination
$$g = r_1g^{\alpha_1}_0+t_1g^{\alpha_1}_1+
\dots +
r_Mg^{\alpha_M}_0+t_Mg^{\alpha_M}_1.$$
over $\bbZ$, for some $\alpha_1<\dots<\alpha_M<\aleph_1$.
Consider the partial sum $r_1g^{\alpha_1}_0+t_1g^{\alpha_1}_1$.
The minimal solution function belongs to $M_{\alpha_1}$,
and using the fact that $h_{\alpha_1}$ increases much faster
than members of $M_{\alpha_1}$, we find infinitely many $j$ such that
$$|r_1g^{\alpha_1}_0+t_1g^{\alpha_1}_1|\rest[0,j)\le j.$$
By induction, we assume that for infinitely many $j$,
the $m-1$-st partial sum satisfies
\begin{equation}\label{goodman}
|r_1g^{\alpha_1}_0+t_1g^{\alpha_1}_1+ \dots + r_{m-1}g^{\alpha_{m-1}}_0+t_{m-1}g^{\alpha_{m-1}}|\rest[0,j)\le cj,
\end{equation}
where $c$ is some constant, and prove the same assertion for $m$-th partial sum.
The set of $j$-s which satisfy \eqref{goodman} defines a function which belongs to $M_{\alpha_{m-1}}$,
and consequently almost each interval $\intvl{h_{\alpha_{m}}}{n}$ contains such a $j$.
Take $n$ such that $(a_n,b_n)=(r_m,t_m)$, and take $j$ satisfying \eqref{goodman}.
On the interval $[h_{\alpha_{m}}(n),j)$, $r_{m}g^{\alpha_{m}}_0+t_{m}g^{\alpha_{m}}_1$ is $0$, and thus the
$m$-th partial sum is is equal to the $m-1$-st partial sum, and the same bound $cj$ applies on that interval.
To take care of $[0,h_{\alpha_{m}}(n))$, we modify the above argument so that $h_{\alpha_{m}}(n+1)\le j$,
and use the fact that $h_{\alpha_{m}}(n)$ is much smaller than $h_{\alpha_{m}}(n+1)$ in a direct calculation.

At the end, there will be infinitely many $j$ such that the $M$-th partial sum, which is equal to $g$,
will be bounded on $[0,j)$ by some constant multiple of $j$, which is bounded, for example, by $j^2$.

\bpf[Proof of Theorem \ref{mini}]
Fix a partition of $\N$ into infinitely many infinite sets $I_l$, $l\in\N$.
Replacing each $I_l$ with the set $\{2n,2n+1 : n\in I_l\}$, we may assume that
for each even $n$, $n\in I_l$ if, and only if, $n+1\in I_l$.
Enumerate $\bbZ^{2}$ as $\sseq{(a_n,b_n)}$, such that the sequence $\{(a_n,b_n)\}_{n\in I_l}$
is constant for each $l$.
Fix an enumeration $\{d_\alpha : \alpha<\aleph_1\}$ of all increasing members of $\NN$.

We carry out a construction by induction on $\alpha<\aleph_1$.
Step $\alpha$:

For each $m$, take a solution to the homogeneous linear equation $a_mx+b_my=0$ over $\bbQ$.
Multiplying $(x,y)$ by a large enough integer multiple of the common denominator of $x$ and $y$,
we may assume that $x,y\in\bbZ$ and $\max\{|x|,|y|\}\ge N$ for any prescribed $N$.
Using that, define nondecreasing functions $\vphi_{\alpha,m}\in\NN$, $m\in\N$, by
$$\vphi_{\alpha,m}(n)=\min\left\{\max\{|x|,|y|\} : \begin{array}{l}
x,y\in\bbZ,\  a_mx+b_my=0,\\
\max\{|x|,|y|\}\ge d_\alpha(n)
\end{array}
\right\},$$
and consequently define $\vphi_\alpha\in\NN$ by
$$\vphi_\alpha(n)=\max\{\vphi_{\alpha,m}(n) : m\le n\}.$$

Let $M_\alpha\sbst\Bgp$
be the smallest set (with respect to inclusion) containing $\vphi_\alpha$ and all
functions defined in stages $<\alpha$, and such that $M_\alpha$ is closed under all operations
relevant for the proof.
For example, closing $M_\alpha$ under the following operations
suffices:
\be
\itm[(a)]\label{hat} $g(n)\mapsto \hat g(n)=\max\{|g(m)| : m\le n\}$;
\itm[(b)]\label{lt} For each $c\in\N$:
$g(n)\mapsto g_{<c}(n)=\min\{j : n\le j,\ g(j)<c\cdot(j+1)\}$,
whenever $g_{<c}$ is well-defined;
\itm[(c)] $(f,g)\mapsto f+g$;
\itm[(d)] $g\mapsto -g$.
\ee
((c)+(d) mean that $M_\alpha\le\Bgp$.)

By induction, $M_\alpha$ is countable.
Take an increasing $h_\alpha\in\NN$ such that
for each $f\in M_\alpha$,
$$f(h_\alpha(n))<h_\alpha(n+1)$$
for all but finitely many $n$.\footnote{\label{aleph1b}To achieve that,
enumerate $M_\alpha\cap\NN=\sseq{f_n}$,
define $h_\alpha(0)=0$, and inductively for each $n>0$,
define $h_\alpha(n+1)=\max\{h_\alpha(n),f_0(h_\alpha(n)),\dots,f_n(h_\alpha(n))\}$+1.}

Define $g^\alpha_0,g^\alpha_1\in\Bgp$ as follows:
For each $n$, choose $c,d\in \bbZ^{2}$ as in the definition of $\vphi_{\alpha,n}(h_\alpha(n+1))$,
and define $(g^\alpha_0(h_\alpha(n)),g^\alpha_1(h_\alpha(n))) = (c,d)$, so that for all $n$,
\begin{eqnarray}
a_ng^\alpha_0(h_\alpha(n))+b_n g^\alpha_1(h_\alpha(n)) & = & 0\label{8c'};\mbox{ and}\\\nonumber
\max\{|g^\alpha_0(h_\alpha(n))|,|g^\alpha_1(h_\alpha(n))|\}  & = &\vphi_{\alpha,n}(h_\alpha(n+1))
\ge d_\alpha(h_\alpha(n+1)).
\end{eqnarray}
The remaining values of the functions $g^\alpha_i$ are defined by declaring these functions constant on
each interval $\intvl{h_\alpha}{n}$.

Take the generated subgroup $G=\langle g^\alpha_0,g^\alpha_1 : \alpha<\aleph_1\rangle$
of $\Bgp$.
We will show that $G$ is as
required in the theorem.

\subsubsection*{$G^2$ is not Menger-bounded}
Let $f\in\NN$.
Take $\alpha<\aleph_1$ such that $f(n)<d_\alpha(n)$ for all $n$.
For each $m$, let
$n$ be such that $m-1\in\intvl{h_\alpha}{n}$. As each function $g^\alpha_i$ is
constant on the interval $\intvl{h_\alpha}{n}$, we have by \eqref{8c'} that
\begin{eqnarray*}
\lefteqn{\max\{|g^\alpha_0(m-1)|,|g^\alpha_1(m-1)|\} =}\\
& = & \max\{|g^\alpha_0(h_\alpha(n))|,|g^\alpha_1(h_\alpha(n))|\} \ge d_\alpha(h_\alpha(n+1))\ge d_\alpha(m)> f(m).
\end{eqnarray*}
This violates Theorem \ref{kMenBdd}(4)
for $k=2$.

\subsubsection*{$G$ is Menger-bounded}
Take $f(n)=n^2$. We will prove that $f$ is as
required in Theorem \ref{kMenBdd}(4).

Fix $g\in G$.
Then there are $M\in\N$, $\alpha_1<\dots<\alpha_M<\aleph_1$,
and integers $r_1,t_1,\dots,r_M,t_M$, such that
$$g = r_1g^{\alpha_1}_0+t_1g^{\alpha_1}_1+
\dots +
r_Mg^{\alpha_M}_0+t_Mg^{\alpha_M}_1.$$
Let $g_0=0$, and for each $m=1,\dots,M$ define
\beq{12'}
g_{m} =
r_1g^{\alpha_1}_0+t_1g^{\alpha_1}_1+
\dots +
r_mg^{\alpha_m}_0+t_mg^{\alpha_m}_1.
\eeq
We prove, by induction on $m=0,\dots,M$, that for an appropriate constant $c_m$,
we have (using the notation in (a) on page \pageref{hat}) that
$$\hat g_{m}(j)\le c_m\cdot(j+1)$$
for infinitely many $j$.

The case $m=0$ is trivial. We show how to move from $m-1$ to $m$.
Assume that
$$J_{m-1} = \{j : \hat g_{m-1}(j) \le c_{m-1}\cdot (j+1)\},$$
is infinite.

By \eqref{8c'}, for each $n>0$,
\begin{eqnarray*}
\lefteqn{\max\{|g^{\alpha_m}_0(h_{\alpha_m}(n-1))|,|g^{\alpha_m}_1(h_{\alpha_m}(n-1))|\}=}\\
& = & \vphi_{{\alpha_m},n-1}(h_{\alpha_m}(n))\le \vphi_{\alpha_m}(h_{\alpha_m}(n)).
\end{eqnarray*}
As $\vphi_{\alpha_m}$ and $h_{\alpha_m}$ are nondecreasing,
$$
\max\{\hat g^{\alpha_m}_0(h_{\alpha_m}(n-1)),\hat g^{\alpha_m}_1(h_{\alpha_m}(n-1))\}\le
\vphi_{\alpha_m}(h_{\alpha_m}(n)),
$$
and since $\vphi_{\alpha_m}\in M_{\alpha_m}$,
\begin{eqnarray}
\lefteqn{\max\{\hat g^{\alpha_m}_0(h_{\alpha_m}(n-1)),\hat g^{\alpha_m}_1(h_{\alpha_m}(n-1))\}\le}\label{14'}\\\nonumber
& \le & \vphi_{\alpha_m}(h_{\alpha_m}(n))<h_{\alpha_m}(n+1).
\end{eqnarray}
for all but finitely many $n$.

As $\alpha_1,\dots,\alpha_{m-1}<\alpha_m$ and $M_{\alpha_m}\le\Bgp$,
$g_{m-1}\in M_{\alpha_m}$. Thus, $\tilde g=\hat g_{m-1}\in M_{\alpha_m}$.
As $J_{m-1}$ is infinite, we have (using the notation of (b) on page \ref{lt})
that
$\tilde g_{<c_{m-1}}(n) = \min\{j : n\le j\in J_{m-1}\}$
is well defined, and $\tilde g_{<c_{m-1}}\in M_{\alpha_m}$.
Consequently,
$\tilde g_{<c_{m-1}}(h_{\alpha_m}(n+1))<h_{\alpha_m}(n+2)$
for all but finitely many $n$. In other words,
for each large enough $n$, there is $j\in J_{m-1}$ such that
\beq{j}
h_{\alpha_m}(n+1)\le j<h_{\alpha_m}(n+2).
\eeq

Let $l$ be such that for each $n\in I_l$, $(a_n,b_n)=(r_m,t_m)$.
For each large enough even $n\in I_l$:
$(a_n,b_n)=(a_{n+1},b_{n+1})=(r_m,t_m)$, and thus by \eqref{8c'},
\begin{eqnarray*}
r_m g^{\alpha_m}_0(h_{\alpha_m}(n))+t_m g^{\alpha_m}_1(h_{\alpha_m}(n)) & = & 0;\\
r_mg^{\alpha_m}_0(h_{\alpha_m}(n+1))+t_mg^{\alpha_m}_1(h_{\alpha_m}(n+1)) & = & 0.
\end{eqnarray*}
By \eqref{12'},
\beq{eq2'}
g_{m}\rest[h_{\alpha_m}(n),h_{\alpha_m}(n+2)) = g_{m-1}\rest[h_{\alpha_m}(n),h_{\alpha_m}(n+2)).
\eeq
Fix $j$ as in \eqref{j}.
Let $p\in[0,j+1)$.

\subsubsection*{Case 1: $p\ge h_{\alpha_m}(n)$.}
As $j<h_{\alpha_m}(n+2)$,
$$[h_{\alpha_m}(n),j+1)\sbst [h_{\alpha_m}(n),h_{\alpha_m}(n+2)),$$
and by \eqref{eq2'} and the membership $j\in J_{m-1}$,
\beq{14.5'}
|g_{m}(p)|=|g_{m-1}(p)|\le\hat g_{m-1}(j)\le c_{m-1}\cdot(j+1).
\eeq

\subsubsection*{Case 2: $p < h_{\alpha_m}(n)$.}
By the definition of $g_{m}$,
$$|g_{m}(p)| \le |g_{m-1}(p)|+2\max\{|r_m|,|t_m|\}\cdot\max\{|g^{\alpha_m}_0(p)|,|g^{\alpha_m}_1(p)|\}.$$
As $p<h_{\alpha_m}(n)\le j\in J_{m-1}$,
$|g_{m-1}(p)|\le \hat g_{m-1}(j)\le c_{m-1}\cdot (j+1)$.
Us\-ing $p\le h_{\alpha_m}(n)-1$, \eqref{14'}, and $h_{\alpha_m}(n+1)\le j$,
we obtain
\begin{equation*}
|g^{\alpha_m}_i(p)|\le \hat g^{\alpha_m}_i(h_{\alpha_m}(n)-1)=
\hat g^{\alpha_m}_i(h_{\alpha_m}(n-1))< h_{\alpha_m}(n+1)\le j.
\end{equation*}
for each $i=0,1$. Together with \eqref{14.5'}, we have that
\begin{eqnarray*}
|g_{m}(p)|
& \le & |g_{m-1}(p)|+2\max\{|r_m|,|t_m|\}\cdot\max\{|g^{\alpha_m}_0(p)|,|g^{\alpha_m}_1(p)|\}\\
& \le & c_{m-1}\cdot(j+1)+2\max\{|r_m|,|t_m|\} j\\
& \le & c_{m-1}\cdot(j+1)+2\max\{|r_m|,|t_m|\}\cdot (j+1)\\
& = & (c_{m-1}+2\max\{|r_m|,|t_m|\})\cdot (j+1).
\end{eqnarray*}
Take $c_m = c_{m-1}+2\max\{|r_m|,|t_m|\}$.

We have proved that for almost all even $n\in I_l$, there is $j\in[h_{\alpha_m}(n+1),h_{\alpha_m}(n+2))$
such that $j\in J_m$. There are infinitely many even $n\in I_l$, and therefore $J_m$ is infinite.
This completes the inductive proof.

Now, for each $j$ in the infinite set $J_M$ such that $c_M\le j$,
$$|g|\rest[0,j+1)\le \hat g(j)\le c_M\cdot(j+1)\le (j+1)^2=f(j+1).$$
By Theorem \ref{kMenBdd}, $G$ is Menger-bounded.
\epf

It is rather straightforward to extend the above proof
to get for each $k$, a group $G\le\Bgp$ such that $G^k$ is
Menger-bounded, but $G^{k+1}$ is not.
To see that, have a quick look at the proof of Theorem \ref{main}.

\section{The main theorem}

Our main Theorem \ref{main} requires a weak portion of \CH{},
that is best stated in terms of cardinal characteristics of the
continuum.
An excellent introduction to the topic is \cite{BlassHBK}.
However, we give a self-contained treatment.

For $f,g\in\NN$, $f\le^* g$ means: $f(n)\le g(n)$ for all but
finitely many $n$. A subset $Y$ of $\NN$ is \emph{bounded} if
there is $g\in\NN$ such that $f\le^* g$ for all $f\in Y$.
At the other extreme, a subset $Y$ of $\NN$ is \emph{dominating} if
for each $f\in\NN$ there is $g\in Y$ such that $f\le^* g$.

$\fb$ is the minimal cardinality of an unbounded subset of $\NN$,
and $\fd$ is the minimal cardinality of a dominating subset of $\NN$.
An argument as in Footnote \ref{aleph1b} shows that $\aleph_1\le\fb$.
Thus, $\aleph_1\le\fb\le\fd\le 2^{\alephes}$.
The hypothesis $\fb=\fd$ is strictly weaker than
\CH{} \cite{BlassHBK}. By inspection, one can see that for the proof of Theorem
\ref{mini}, it suffices to assume that $\fb=\fd$.
To extend this observation further, we introduce the following
new cardinal characteristics.

\bdfn\label{dprime}
Fix a partition $\cP=\{I_l : l\in\N\}$ of $\N$ such that for each
$l$, there are infinitely many $n$ such that $n,n+1\in I_l$.
For $f\in\NN$ and an increasing $h\in\NN$, write
$$[f\ll h] = \{n : f(h(n))<h(n+1)\}.$$
$\fd'(\cP)$ is the cardinal such that the following are equivalent:
\be
\itm $\kappa<\fd'(\cP)$;
\itm For each $Y\sbst\NN$ such that $|Y|=\kappa$,
there is an increasing $h\in\NN$ such that for each $f\in Y$,
$$(\forall l)(\exists^\oo n)\ n,n+1\in I_l\cap [f\ll h].$$
\ee
\edfn

Clearly, $\fb\le\fd'(\cP)\le\fd$ for each $\cP$.
We first point out that the hypothesis ``there is $\cP$ such that $\fd'(\cP)=\fd$''
is strictly weaker than the hypothesis $\fb=\fd$.
Let $\cov(\cM)$ be the minimal cardinality of a cover of $\NN$ by
meager (first category) sets. It is consistent that $\fb<\cov(\cM)=\fd$ \cite{BlassHBK}.

\blem
For each $\cP$, $\cov(\cM)\le\fd'(\cP)$.
\elem
\bpf
Fix a partition $\cP=\{I_l : l\in\N\}$ of $\N$ such that for each
$l$, there are infinitely many $n$ such that $n,n+1\in I_l$.
Let $\NNup$ be the set of all increasing elements of $\NN$.
$\NNup$ is homeomorphic to $\NN$. It therefore suffices to find
a cover of $\NNup$ by $\fd'(\cP)$ many nowhere-dense subsets of $\NNup$.

Take $Y\sbst\NN$ such that $|Y|=\fd'(\cP)$, and
such that Definition \ref{dprime}(2)
fails for $Y$, that is: For each $h\in\NNup$, there are $f\in Y$ and $l$ such that
$$(\forall^\oo n)\ n,n+1\in I_l\to f(h(n))\ge h(n+1)\mbox{ or } f(h(n+1))\ge h(n+2).$$
For $f\in Y$ and $l,m\in\N$, let
$$Y_{f,l,m} = \left\{h\in\NNup :
\begin{array}{l}
(\forall n\ge m)\ n,n+1\in I_l\to\\
f(h(n))\ge h(n+1)\mbox{ or } f(h(n+1))\ge h(n+2)
\end{array}\right\}.$$
$Y_{f,l,m}$ is nowhere dense in $\NNup$: Given $k$ and an increasing finite sequence $s\in\N^k$, let
$n\ge\max\{k,m\}$ be such that $n,n+1\in I_l$. Let $\tilde s$ be an extension of $s$ to an increasing
sequence of length $n+3$, such that $f(\tilde s(n))<\tilde s(n+1)$ and $f(\tilde s(n+1))<\tilde s(n+2)$.
Then $Y_{f,l,m}\cap\bo{\tilde s}=\emptyset$.

As $\Union\{Y_{f,l,m} : f\in Y, l,m\in\N\}=\NNup$, $\cov(\cM)\le\fd'(\cP)\cdot\alephes=\fd'(\cP)$.
\epf

A more thorough analysis of the cardinals
$\fd'(\cP)$ is carried out by Mildenberger \cite{MilBG}.

\bthm\label{main}
Assume that there is $\cP$ such that $\fd'(\cP)=\fd$.
Then for each $k$, there is a group $G\le\Bgp$ such that
$G^k$ is Menger-bounded, but $G^{k+1}$ is not Menger-bounded.
\ethm
\bpf
Fix a partition $\cP=\{I_l : l\in\N\}$ of $\N$ such that for each
$l$, there are infinitely many $n$ such that $n,n+1\in I_l$,
and such that $\fd'(\cP)=\fd$.

Enumerate $\bbZ^{k\x (k+1)}$ as $\sseq{A_n}$, such that the sequence $\{A_n\}_{n\in I_l}$
is constant for each $l$.
Fix a dominating family of increasing functions $\{d_\alpha : \alpha<\fd\}\sbst\NN$.
For $v=(v_0,\dots,v_k)\in \bbZ^{k+1}$, write $\|v\|$ or $\|v_0,\dots,v_k\|$ for $\max\{|v_0|,\dots,|v_k|\}$
(the supremum norm of $v$).

We carry out a construction by induction on $\alpha<\fd$.
Step $\alpha$:
Define functions $\vphi_{\alpha,m}\in\NN$, $m\in\N$, by
$$\vphi_{\alpha,m}(n)=\min\{\|v\| : v\in \bbZ^{k+1},\ \|v\|\ge d_\alpha(n),\ A_mv=\vzero\}.$$
Also, define $\vphi_\alpha\in\NN$ by
\beq{2}
\vphi_\alpha(n)=\max\{\vphi_{\alpha,m}(n) : m\le n\}.
\eeq

Let $M_\alpha\sbst\Bgp$
be the smallest set (with respect to inclusion) containing $\vphi_\alpha$ and all
functions defined in stages $<\alpha$, and such that $M_\alpha$ is closed under all operations
relevant for the proof.
For example, closing $M_\alpha$ under the following operations
suffices:
\be
\itm[(a)] $g(n)\mapsto \hat g(n)=\max\{|g(m)| : m\le n\}$;
\itm[(b)] $(g(n),f(n))\mapsto\max\{|g(n)|,|f(n)|\}$;
\itm[(c)] For each $c\in\N$:
$g(n)\mapsto g_{<c}(n)=\min\{j : n\le j,\ g(j)<c\cdot(j+1)\}$,
whenever $g_{<c}$ is well-defined.
\itm[(d)] $(f,g)\mapsto f+g$;
\itm[(e)] $g\mapsto -g$.
\ee
There are countably many such operations, and by induction,
$|M_\alpha|\le\max\{\alephes,|\alpha|\}<\fd=\fd'(\cP)$.
By the definition of $\fd'(\cP)$, there is
an increasing $h_\alpha\in\NN$ such that
for each $f\in M_\alpha\cap\NN$,
\beq{7}
(\forall l)(\exists^\oo n)\ n,n+1\in I_l\cap[f\ll h_\alpha].
\eeq
Define $k+1$ elements $g^\alpha_0,\dots,g^\alpha_k\in\Bgp$ as follows:
For each $n$, let $v\in \bbZ^{k+1}$ be a witness for the definition of $\vphi_{\alpha,n}(h_\alpha(n+1))$,
namely,
\begin{eqnarray}
\vphi_{\alpha,n}(h_\alpha(n+1)) & = & \|v\|\ge d_\alpha(h_\alpha(n+1))\label{8a}\\
A_nv & = & \vzero.\label{8b}
\end{eqnarray}
and define
\beq{8c}
\mx{g^\alpha_0(h_\alpha(n))\\
\vdots\\
g^\alpha_k(h_\alpha(n))}
=v,$$
so that
$$A_n\cdot \mx{g^\alpha_0(h_\alpha(n))\\
\vdots\\
g^\alpha_k(h_\alpha(n))}
=\vzero.
\eeq
The remaining values of the functions $g^\alpha_i$ are defined by declaring these functions constant on
each interval $\intvl{h_\alpha}{n}$.
By \eqref{8a} and \eqref{8c},
\beq{eq1}
\|g^\alpha_0(h_\alpha(n)),\dots,\allowbreak g^\alpha_k(h_\alpha(n))\|=\vphi_{\alpha,n}(h_\alpha(n+1)).
\eeq
for all $n$.

Take the generated subgroup $G=\langle g^\alpha_0,\dots,g^\alpha_k : \alpha<\fd\rangle$
of $\Bgp$.
We will show that $G$ is as
required in the theorem.

\subsubsection*{$G^{k+1}$ is not Menger-bounded}
We use Theorem \ref{kMenBdd}.
Let $f\in\NN$.
Take $\alpha<\fd$ such that $f<^* d_\alpha$, and set $F=\{g^\alpha_0,\dots,g^\alpha_k\}\in[G]^{k+1}$.
For each large enough $m$: $f(m)<d_\alpha(m)$. Fix such $m$. Let
$n$ be such that $m-1\in\intvl{h_\alpha}{n}$. As each function $g^\alpha_i$ is
constant on the interval $\intvl{h_\alpha}{n}$, and using \eqref{eq1} and \eqref{8a}, we have that
\begin{eqnarray*}
\lefteqn{\|g^\alpha_0(m-1),\dots,g^\alpha_k(m-1)\|=}\\
& = & \|g^\alpha_0(h_\alpha(n)),\dots,g^\alpha_k(h_\alpha(n))\| = \vphi_{\alpha,n}(h_\alpha(n+1))\ge\\
& \ge & d_\alpha(h_\alpha(n+1))\ge d_\alpha(m)> f(m).
\end{eqnarray*}
This violates Theorem \ref{kMenBdd}(4)
for the power $k+1$.

\subsubsection*{$G^{k}$ is Menger-bounded}
Take $f(n)=n^2$. Clearly, $f$ dominates all functions
$f_c(n)=c\cdot n$, $c\in\N$. We will prove that $f$ is as
required in Theorem \ref{kMenBdd}(4).

Fix $F=\{g_0,\dots,g_{k-1}\}\sbst G$.
Then there are $M\in\N$, $\alpha_1<\dots<\alpha_M<\fd$,
and matrices $B_1,\dots,B_M\in\bbZ^{k\x (k+1)}$, such that
$$\mx{g_0\\\vdots\\ g_{k-1}} =
B_1\mx{g^{\alpha_1}_0\\\vdots\\g^{\alpha_1}_k}+
\dots +
B_M\mx{g^{\alpha_M}_0\\\vdots\\g^{\alpha_M}_k}.$$
Let $g_{0,0}=\dots=g_{k-1,0}=0$, and for each $m=1,\dots,M$ let
\beq{12}
\mx{g_{0,m}\\\vdots\\ g_{k-1,m}} =
B_1\mx{g^{\alpha_1}_0\\\vdots\\g^{\alpha_1}_k}+
\dots +
B_m\mx{g^{\alpha_m}_0\\\vdots\\g^{\alpha_m}_k}.
\eeq
We prove, by induction on $m=0,\dots,M$, that for an appropriate constant $c_m$,
there are infinitely many $j$ such that
$$\|\hat g_{0,m}(j), \dots, \hat g_{k-1,m}(j)\|\le c_m\cdot(j+1).$$
By the definition of $f$, this suffices.

The case $m=0$ is trivial. We show how to move from $m-1$ to $m$.
Assume that
$$J_{m-1} = \{j : \|\hat g_{0,m-1}(j), \dots, \hat g_{k-1,m-1}(j)\| \le c_{m-1}\cdot (j+1)\},$$
is infinite.

As $\alpha_1,\dots,\alpha_{m-1}<\alpha_m$, we have by \eqref{12} that
$g_{0,m-1},\dots,g_{k-1,m-1}\in M_{\alpha_m}$. By (a),(b),(c), the functions
$$g(n)=\|\hat g_{0,m-1}(n), \dots, \hat g_{k-1,m-1}(n)\|$$
and $g_{<c_{m-1}}$ both belong to $M_{\alpha_m}$.
Note that
\beq{13.5}
g_{<c_{m-1}}(n) = \min\{j : n\le j\in J_{m-1}\},
\eeq
and is therefore well defined.
Thus, $\max\{g_{<c_{m-1}},\vphi_{\alpha_m}\}\in M_{\alpha_m}$.

For each $i\le k$ and each $n>0$, as $n-1\le h_{\alpha_m}(n)$, we have by \eqref{eq1} that
$$|g^{\alpha_m}_i(h_{\alpha_m}(n-1))|\le \vphi_{{\alpha_m},n-1}(h_{\alpha_m}(n))\le \vphi_{\alpha_m}(h_{\alpha_m}(n)).$$
As $\vphi_{\alpha_m}$ and $h_{\alpha_m}$ are nondecreasing,
\beq{14}
\|\hat g^{\alpha_m}_0(h_{\alpha_m}(n-1)),\dots,\hat g^{\alpha_m}_k(h_{\alpha_m}(n-1))\|\le
\vphi_{\alpha_m}(h_{\alpha_m}(n)).
\eeq
Thus, if $l$ is such that for each $n\in I_l$, $A_n=B_m$, we have by \eqref{13.5} and \eqref{14} that
\begin{eqnarray*}
I & = &
\left\{n :
\begin{array}{l}
A_n=A_{n+1}=B_m,\\
(\exists j\in J_{m-1})\ h_{\alpha_m}(n+1)\le j<h_{\alpha_m}(n+2),\\
\|\hat g^{\alpha_m}_0(h_{\alpha_m}(n-1)),\dots,\hat g^{\alpha_m}_k(h_{\alpha_m}(n-1))\|< h_{\alpha_m}(n+1)
\end{array}
\right\}\\
& \spst &
\left\{n :
\begin{array}{l}
n,n+1\in I_l,\\
g_{<c_{m-1}}(h_{\alpha_m}(n+1))<h_{\alpha_m}(n+2),\\
\vphi_{\alpha_m}(h_{\alpha_m}(n))<h_{\alpha_m}(n+1)
\end{array}
\right\}\\
& \spst &
\{n : n,n+1\in I_l\cap[\max\{g_{<c_{m-1}},\vphi_{\alpha_m}\}\ll h_{\alpha_m}]\}.
\end{eqnarray*}
As $\max\{g_{<c_{m-1}},\vphi_{\alpha_m}\}\in M_{\alpha_m}$,
we have by the definition of $h_{\alpha_m}$ \eqref{7} that
the last set is infinite, and therefore so is $I$.

Let $n\in I$. Then $A_n = A_{n+1} = B_m$, and thus by \eqref{8b} and \eqref{8c},
$$B_m\cdot\mx{g^{\alpha_m}_0(h_{\alpha_m}(n))\\\vdots\\g^{\alpha_m}_k(h_{\alpha_m}(n))}=
B_m\cdot
\mx{g^{\alpha_m}_0(h_{\alpha_m}(n+1))\\\vdots\\g^{\alpha_m}_k(h_{\alpha_m}(n+1))}=\vzero.$$
By \eqref{12}, for each $i<k$,
\beq{eq2}
g_{i,m}\rest[h_{\alpha_m}(n),h_{\alpha_m}(n+2)) = g_{i,m-1}\rest[h_{\alpha_m}(n),h_{\alpha_m}(n+2)).
\eeq
As $n\in I$, there is $j\in J_{m-1}$ such that $h_{\alpha_m}(n+1)\le j<h_{\alpha_m}(n+2)$,
and
\beq{eq3}
\|\hat g^{\alpha_m}_0(h_{\alpha_m}(n-1)),\dots,\hat g^{\alpha_m}_k(h_{\alpha_m}(n-1))\|< h_{\alpha_m}(n+1).
\eeq
Let $p\in[0,j+1)$.

\subsubsection*{Case 1: $p\ge h_{\alpha_m}(n)$.}
As $j<h_{\alpha_m}(n+2)$,
$$[h_{\alpha_m}(n),j+1)\sbst [h_{\alpha_m}(n),h_{\alpha_m}(n+2)),$$
and by \eqref{eq2} and the membership $j\in J_{m-1}$,
\beq{14.5}
|g_{i,m}(p)|=|g_{i,m-1}(p)|\le\hat g_{i,m-1}(j)\le c_{m-1}\cdot(j+1)
\eeq
for all $i<k$.

\subsubsection*{Case 2: $p < h_{\alpha_m}(n)$.}
Let $C$ be the maximal absolute value of a coordinate of $B_m$.
For all $i<k$, by the definition of $g_{i,m}$,
\beq{15}
|g_{i,m}(p)| \le |g_{i,m-1}(p)|+(k+1)C\cdot\max\{|g^{\alpha_m}_i(p)| : i\le k\}.
\eeq
As $p<h_{\alpha_m}(n)\le j\in J_{m-1}$,
$|g_{i,m-1}(p)|\le \hat g_{i,m-1}(j)\le c_{m-1}\cdot (j+1)$.
Using $p\le h_{\alpha_m}(n)-1$,\eqref{eq3}, $g^{\alpha_m}_i$ being constant on
$[h_{\alpha_m}(n-1),h_{\alpha_m}(n))$, and $h_{\alpha_m}(n+1)\le j$,
we obtain
\begin{equation*}
|g^{\alpha_m}_i(p)|\le \hat g^{\alpha_m}_i(h_{\alpha_m}(n)-1)=
\hat g^{\alpha_m}_i(h_{\alpha_m}(n-1))< h_{\alpha_m}(n+1)\le j.
\end{equation*}
for each $i\le k$. Together with \eqref{14.5}, we have that
\begin{eqnarray*}
|g_{i,m}(p)|
& \le & |g_{i,m-1}(p)|+(k+1)C\cdot\max\{|g^{\alpha_m}_i(p)| : i\le k\}\\
& \le & c_{m-1}\cdot(j+1)+(k+1)C j\\
& \le & c_{m-1}\cdot(j+1)+(k+1)C\cdot (j+1)\\
& = & (c_{m-1}+(k+1)C)\cdot (j+1).
\end{eqnarray*}
Take $c_m = c_{m-1}+(k+1)C$.
We have proved that for each $n\in I$ there is $j\in[h_{\alpha_m}(n+1),h_{\alpha_m}(n+2))$
such that $j\in J_m$. $I$ is infinite, and therefore so is $J_m$.
This completes the inductive proof, and consequently the proof of Theorem \ref{main}.
\epf

\brem
Mildenberger has recently proved that our assumption in Theorem \ref{main}
can be weakened to $\fd\le\fr$  \cite{MilBG}.
\erem

\subsection*{Acknowledgements}
We thank Heike Mildenberger and Lyubomyr Zdomskyy for their useful comments.
Following the suggestions of the referee, we have
substantially expanded the introductory and motivational
parts of this paper.

\Newpage

\ed